\documentclass[12pt]{article}

\usepackage[english]{babel}
\usepackage{amssymb,a4wide}
\usepackage{doublespace}
\usepackage{amsmath,xypic,epsf,mathrsfs,theorem}

\usepackage[all]{xy}
\newdir{ (}{{}*!/-5pt/\dir^{(}}
\newdir{ )}{{}*!/-5pt/\dir_{(}}

\newtheorem{theorem}{Theorem}[subsection]
\newtheorem{lemma}[theorem]{Lemma}
\newtheorem{corollary}[theorem]{Corollary}

\newtheorem{proposition}[theorem]{Proposition}

{\theorembodyfont{\upshape}
\newtheorem{definition}[theorem]{Definition}

}

\numberwithin{equation}{section}
\numberwithin{theorem}{section}

\newcommand{\image}{\mathrm{Im}}

\newcommand{\diag}{\mathrm{diag}}

\newcommand{\GL}{\mathrm{GL}}

\newcommand{\ev}{\mathrm{ev}}
\newcommand{\Aff}{\mathrm{Aff}}

\newcommand{\RR}{{\rm RR}}
\newcommand{\tr}{{\rm tr}}
\newcommand{\cM}{{\cal M}}
\newcommand{\ep}{{\varepsilon}}

\newcommand{\cU}{{\cal U}}
\newcommand{\cP}{{\cal P}}

\newcommand{\cO}{{\cal O}}

\newcommand{\cZ}{{\cal Z}}
\newcommand{\C}{{\mathbb C}}
\newcommand{\Z}{{\mathbb Z}}
\newcommand{\N}{{\mathbb N}}

\newcommand{\cK}{{\cal K}}

\newcommand{\R}{{\mathbb R}}
\newcommand{\Cs}{{$C^*$-al\-ge\-bra}}

\newcommand{\sh}{{$^*$-ho\-mo\-mor\-phism}}

\newcommand\eqdef{{\;\overset{\mbox{\scriptsize def}}{=}\;}}

\newenvironment{proof}[1][Proof:]
{\begin{trivlist}\item[]\textbf{#1} }
{\hbox{}\nobreak\hfill\quad\hbox{$\square$}\end{trivlist}}

\begin{document}

\title{The stable and the real rank of $\cZ$-absorbing $C^*$-algebras}

\author{Mikael R\o rdam}
\date{}
\maketitle

\begin{abstract} \noindent Suppose that $A$ is a \Cs{} for which $A
  \cong A \otimes \cZ$, where $\cZ$ is the Jiang--Su algebra: a
  unital, simple, stably finite, separable, nuclear, 
  infinite dimensional \Cs{} with the same
  Elliott invariant as the complex numbers. We show that:
\begin{enumerate}
\item The Cuntz semigroup $W(A)$ of equivalence classes of positive
  elements in matrix algebras over $A$ is almost 
unperforated\footnote{Almost perforation is a natural extension of the 
notion of weak unperforation for \emph{simple} 
ordered abelian (semi-)groups, see Section~3.}.
\item If $A$ is exact, then $A$ is purely infinite if and only if $A$
  is traceless.
\item If $A$ is separable and nuclear, then $A \cong A \otimes
  \cO_\infty$ if and only if $A$ is traceless.
\item If $A$ is simple and unital, then the stable rank of $A$ is one if
  and only if $A$ is finite.
\end{enumerate}
We also characterise when $A$ is of real rank zero.
\end{abstract}

\section{Introduction} \label{sec:intro} 

\noindent Jiang and Su gave in their paper \cite{JiaSu:Z} a
classification of simple inductive limits of direct sums of dimension drop
\Cs s. (A dimension drop \Cs{} is a certain sub-\Cs{} of $M_n(C([0,1]))$,
a precise definition of which is given in the next section.) 
They prove that inside this class there exists a
unital, simple, infinite dimensional \Cs{} $\cZ$ whose Elliott
invariant is isomorphic to the Elliott invariant of the complex
numbers, that is,  
$$(K_0(\cZ), K_0(\cZ)^+, [1]) \cong (\Z, \Z^+, 1), \qquad K_1(\cZ)=0,
\qquad T(\cZ) = \{\tau\},$$
where $\tau$ is the unique tracial state on $\cZ$. They proved that
$\cZ \cong \cZ \otimes \cZ \cong \bigotimes_{n=1}^\infty \cZ$, and
that $A \otimes \cZ \cong A$ if $A$ is a simple, unital, infinite
dimensional AF-algebra or if $A$ is a unital Kirchberg algebra. Toms
and Winter have in a paper currently under preparation extended the
latter result by showing that $A \otimes \cZ \cong A$ for all
approximately divisible \Cs s. Toms and Winter note that upon
combining results from \cite{JiaSu:Z} with \cite[Theorem~8.2]{KirRor:pi2}
one obtains that a separable \Cs{} $A$ is $\cZ$-absorbing if and only
if there is a unital embedding of $\cZ$ into the relative commutant
$\cM(A)_\omega \cap A'$,
where $\cM(A)_\omega$ is the ultrapower, associated with a free filter
$\omega$ on $\N$, of the multiplier algebra, $\cM(A)$, of $A$. This
provides a partial answer to the question raised by Gong, Jiang, and
Su in \cite{GongJiangSu:Z} if one can give an intrinsic description of
which (separable, nuclear) \Cs s absorb $\cZ$. 

Gong, Jiang, and Su prove in
\cite{GongJiangSu:Z} that $(K_0(A),K_0(A)^+) \cong (K_0(A \otimes \cZ),
K_0(A \otimes \cZ)^+)$
if and only if $K_0(A)$ is weakly
unperforated as an ordered group, when $A$ is a simple \Cs{}; and 
hence that $A$ and $A \otimes 
\cZ$ have isomorphic Elliott invariant if $A$ is simple
with weakly unperforated $K_0$-group. This result indicates that $A
\cong A \otimes \cZ$ whenever $A$ is ``classifiable'' in the sense of
Elliott (see Elliott, \cite{Ell:classprob}, or \cite{Ror:encyc} by the
author). 

The results quoted above show on the one hand that surprisingly many
\Cs s, including for example the irrational rotation \Cs s, absorb the
Jiang--Su algebra, but on the other hand that not all simple, unital,
nuclear, separable \Cs s are $\cZ$-absorbing. Villadsen's example from
\cite{Vil:perforation} of a simple, unital AH-algebra whose
$K_0$-group is not weakly unperforated cannot absorb $\cZ$. The
example by the author in \cite{Ror:simple} of a simple, unital,
nuclear, separable \Cs{} with a finite and an infinite projection is
prime (i.e., is not the tensor product of two non type I \Cs s), and
is hence not $\cZ$-absorbing. Toms gave in \cite{Toms:example}
an example of a simple, unital ASH-algebra which is not
$\cZ$-absorbing, but which has weakly unperforated $K_0$-group. The
latter two examples (by the author and by Toms) have the same
Elliott invariant as, but are not isomorphic to, their $\cZ$-absorbing
counterparts; and so they serve as counterexamples to the classification
conjecture of Elliott (as it is formulated in \cite[Section~2.2]{Ror:encyc}).  

It appears plausible that the Elliott conjecture holds for all simple,
unital, nuclear, separable $\cZ$-absorbing \Cs s.

In the present paper we begin by showing that the Cuntz semigroup of
equivalence classes of positive elements in a $\cZ$-absorbing \Cs{} is
almost unperforated (a property that for simple ordered abelian 
(semi-)groups coincides with the weak unperforation property, see 
Section~3). We use this to show that the semigroup $V(A)$ of
Murray--von Neumann equivalence classes of projections in a
$\cZ$-absorbing \Cs{} $A$, and in some cases also $K_0(A)$, is almost
(or weakly) unperforated. We show that the stable rank of $A$ is one
if $A$ is a simple, finite, unital $\cZ$-absorbing \Cs{}, thus
answering in the affirmative a question from \cite{GongJiangSu:Z}.
In the last section we characterise when a simple unital
$\cZ$-absorbing \Cs{} is of real rank zero.

\section{Preliminary facts about the Jiang--Su algebra $\cZ$}

\noindent We establish a couple of results that more or less follow
directly from Jiang and Su's paper
\cite{JiaSu:Z} on their \Cs{} $\cZ$. 

For each triple of natural numbers $n,n_0,n_1$, for which $n_0$ and
$n_1$ divides $n$, the \emph{dimension drop \Cs{}} $I(n_0,n,n_1)$ is
the sub-\Cs{} of $C([0,1],M_n)$ consisting of all functions
$f$ such that $f(0) \in \varphi_0(M_{n_0})$ and $f(1) \in
\varphi_1(M_{n_1})$, where $\varphi_j \colon 
M_{n_j} \to M_n$, $j=0,1$, are fixed
unital \sh s. (The \Cs{} $I(n_0,n,n_1)$ is---up to
$^*$-isomorphism---independent on the choice of the \sh s
$\varphi_j$.) The dimension drop \Cs{}  $I(n_0,n,n_1)$ is said to be  
\emph{prime} if $n_0$ and $n_1$ are relatively prime, or equivalently,
if $I(n_0,n,n_1)$ has no projections other than the two trivial ones: $0$ 
and $1$, cf.\ \cite{JiaSu:Z}. 

The \Cs{} $I(n,nm,m)$ can, and will in this paper, be realized as the
sub-\Cs{} of $C([0,1], M_n \otimes M_m)$ consisting of those functions
$f$ for which $f(0) \in M_n \otimes \C$ and $f(1) \in \C \otimes
M_m$. 

A unital \sh{} $\psi \colon I(n_0,n,n_1) \to \cZ$ will here be said to
be \emph{standard}, if 
\begin{equation} \label{eq:tau_0}
\tau(\psi(f)) = \int_0^1 \tr(f(t)) \, dm(t), \qquad f \in I(n_0,n,n_1),
\end{equation}
where $\tau$ is the unique trace on $\cZ$, and where $\tr$ is the
normalised trace on $M_n$. 

The following theorem is essentially contained in Jiang and Su's paper
(\cite{JiaSu:Z}). 

\begin{theorem}[Jiang--Su] \label{thm:I}
Let $n,n_0,n_1$ be a triple of natural numbers where $n_0$ and $n_1$
divide $n$, and where $n_0$ and $n_1$ are relatively prime.
As above, let $\tau$ denote the unique trace on $\cZ$.
\begin{enumerate}
\item For each faithful tracial state $\tau_0$ on $I(n_0,n,n_1)$ 
there exists a unital embedding\\ $\psi \colon I(n_0,n,n_1) \to \cZ$ such
that $\tau \circ \psi = \tau_0$. In particular, there is a standard
unital embedding of $I(n_0,n,n_1)$ into $\cZ$. 
\item Two unital embeddings $\psi_1, \psi_2 \colon I(n_0,n,n_1) \to \cZ$ 
are approximately unitarily equi\-va\-lent if and only if $\tau \circ
\psi_1 = \tau \circ \psi_2$. In particular, $\psi_1$
and $\psi_2$ are approximately unitarily equivalent if they both are
standard.
\end{enumerate} 
\end{theorem}

\begin{proof}
For brevity, denote the prime dimension drop \Cs{} $I(n_0,n,n_1)$
by $I$. 

Find an increasing sequence $B_1 \subseteq B_2 \subseteq B_3
\subseteq \cdots $ of sub-\Cs s of $\cZ$ such that each $B_k$ is
(isomorphic to) a prime
dimension drop algebra of the form $I(n_0(k),n(k),n_1(k))$, and such
that $\bigcup_{k=1}^\infty B_k$ is dense in $\cZ$. Simplicity of $\cZ$
ensures that $n_0(k)$,
$n(k)$, and $n_1(k)$ all tend to infinity as $k$ tends to infinity.

It is shown in \cite[Lemma~2.3]{JiaSu:Z} that  
$K_0( I)$ and $K_0(B_k)$ are infinite cyclic groups
each generated by the class of the unit in the corresponding algebra, 
and $K_1(I)$ and $K_1(B_k)$ are both trivial. This entails that
$KK(\psi_1) = KK(\psi_2)$ for any pair of \emph{unital} \sh s $\psi_1, \psi_2
\colon I \to B_k$.

In both parts of the proof we shall apply the uniqueness theorem,
\cite[Corollary~5.6]{JiaSu:Z}, in Jiang and Su's paper. For each $x
\in [0,1]$ consider the extremal tracial state $\tau_x$ on a dimension
drop
\Cs{} $I(m_0,m,m_1)$ given by $\tau_x(f) = \tr(f(x))$ (where $\tr$ is
the normalised trace on $M_m$). Each self-adjoint 
element $f$ in a dimension drop \Cs{} $I(m_0,m,m_1)$ gives rise to a
function $\hat{f} \in C_\R([0,1])$ defined by $\hat{f}(x) =
\tau_x(f)$. If $f$ is a self-adjoint element in the center of
$I(m_0,m,m_1)$, then $\hat{f} = f$. Every \sh{} $\psi \colon
I(m_0,m,m_1) \to I(m'_0,m',m'_1)$ between two dimension drop \Cs s
induces a positive linear mapping $\psi_* \colon C_\R([0,1]) \to
C_\R([0,1])$ given by $\psi_*(f) = \widehat{{\psi(f)}}$ (when we
identify $C_\R([0,1])$ with the self-adjoint part of the center of
$I(m_0,m,m_1)$).  Let $h_{D,d} \in
C_\R([0,1])$, $d=1, 2, \dots, D$, be the test functions defined in
\cite[(5.5)]{JiaSu:Z} (and previously considered by Elliott).

(i). Let $\tau_0$ be a faithful trace on $I = I(n_0,n,n_1)$. 
Let $F_1 \subseteq F_2 \subseteq
\cdots$ be an increasing sequence of finite subsets of $I$
with dense union. By a one-sided approximate intertwining argument
(after Elliott, see eg.\ \cite[Theorem~1.10.14]{Lin:amenable}) 
it suffices to find a sequence $1
\le m(1) < m(2) < m(3) < \cdots$ of integers, a sequence $\psi_j
\colon  I \to B_{m(j)}$ of unital \sh s, and unitaries
$u_j \in B_{m(j)}$ such that 
$$\|u_{j+1}^*\psi_j(f)u_{j+1} - \psi_{j+1}(f)\| \le 2^{-j}, \qquad 
|\tau(\psi_j(f)) - \tau_0(f)| \le 1/j, \qquad f \in F_j,$$
for all $j \in \N$. It will then follow that
there exist a \sh{} $\psi 
\colon I \to \cZ$ and unitaries $v_j \in \cZ$ such that
$\|v_j^*\psi_j(f)v_j-\psi(f)\|$ tends to zero as $j$ tends to infinity
for all $f \in I$. This will imply that $\tau \circ \psi = \tau_0$.

For each $j$ choose a natural number $D_j$ such that $\|f(s)-f(t)\|
\le 2^{-j}$ for all $f \in F_j$ and for all $s,t \in [0,1]$ with
$|s-t| \le 6/D_j$. Let $G_j$ be the finite set that contains
$F_j$ and the test functions $h_{D_j,d}$, $d=1, \dots, D_j$. Put
$$c_j = \tfrac{2}{3} \min \big\{ \tau_0\big(h_{D_j,d-1}-h_{D_j,d}\big)
\mid d = 2,3, \dots, D_j\big\} >0.$$

By \cite[Corollary 4.4]{JiaSu:Z} --- if $m(j)$ are chosen large enough
--- there exists for each $j$ a unital \sh{} $\psi_j \colon I \to
B_{m(j)}$ such that $|\tau'(\psi_j(f)) - \tau_0(f)| <
\min\{1/j,c_j/2,c_{j-1}/2\}$ for all tracial states $\tau'$ on $B_{m(j)}$
and for all $f \in G_j$. In particular, $|\tau(\psi_j(f))-\tau_0(f)|< 1/j$ 
for $f \in F_j$, and 
\begin{eqnarray*}
&(\psi_j)_*(h_{D_j,d-1}-h_{D_j,d}) \ge c_j, \qquad
(\psi_{j+1})_*(h_{D_j,d-1}-h_{D_j,d}) \ge c_j, &\\
&\| (\psi_{j+1})_*(h_{D_j,d}) - (\psi_{j})_*(h_{D_j,d})
\|_\infty  < c_j,
\end{eqnarray*}
for $d=(1),2,3, \dots, D_j$. It now follows from \cite[Corollary~5.6]{JiaSu:Z}
that there exists a unitary
$u_{j+1}$ in $B_{m(j+1)}$ such that $\|u_{j+1}^*\psi_j(f)u_{j+1} - 
\psi_{j+1}(f)\| \le 2^{-j}$ .

(ii). The ``only if'' part is trivial. Assume that $\tau \circ
\psi_1 = \tau \circ \psi_2$. Take a finite subset $F$ of $I$ and let
$\ep >0$. 

It is shown in
\cite{JiaSu:Z} (and in \cite{EilLorPed:semiproj}) that the dimension drop
\Cs{} $I=I(n_0,n,n_1)$ is semiprojective. We can therefore, 
for some large enough $k_0$, find unital \sh s $\psi_1^{(k)}$, 
$\psi_2^{(k)} \colon I \to B_k$ for each $k \ge k_0$ such that 
\begin{equation} \label{eq1}
\lim_{k \to \infty} \| \psi_j(f) - \psi_j^{(k)}(f)\| = 0, \qquad f \in
I, \; \; j=1,2.
\end{equation}

We assert that
\begin{equation} \label{eq2}
\lim_{k \to \infty} \|(\psi_j^{(k)})_*(h) -
(\tau \circ \psi_j)(h) {\bf 1}\|_\infty = 0, \qquad j=1,2, \quad h
\in C_\R([0,1]),
\end{equation}
when we identify $C_\R([0,1])$ with the self-adjoint portion of the
center of $I = I(n_0,n,n_1)$. Indeed, because $\tau$ is the unique
trace on $\cZ$, the quantity
$$\sup_{\tau' \in T(B_k)} |\tau'(b)-\tau(b)|, \qquad b \in B_\ell,$$
tends to zero as $k$ tends to infinity (with $k \ge \ell$). 
Hence, if we let $\iota_{k,\ell}$ denote the
inclusion mapping $B_\ell \to B_k$, then $\|(\iota_{k,\ell})_*(h) -
\tau(h) {\bf 1}\|_\infty$ tends to zero as $k$ tends to
infinity ($k \ge \ell$). It follows from \eqref{eq1} that
$\|(\iota_{k,\ell} \circ \psi_j^{(\ell)})_*(h) -
(\psi_j^{(k)})_*(h)\|_\infty$ is small if $\ell$ is large (and $k \ge
\ell$). The claim in \eqref{eq2} follows from these facts 
and the identity $(\iota_{k,\ell} \circ
\psi_j^{(\ell)})_* = (\iota_{k,\ell})_* \circ (\psi_j^{(\ell)})_*$.

Choose an
integer $D$ such that $\|f(s)-f(t)\| < \ep/9$ for all $f \in F$ and
for all $s,t \in [0,1]$ with $|s-t| \le 6/D$. 
Each $\psi_j(h_{D,d-1}-h_{D,d})$ is a non-zero and
positive element in $\cZ$, and we can therefore find $c >0$ such that 
$(\tau\circ\psi_j)(h_{D,d-1}-h_{D,d}) \ge 2c$ for $d=2,3, \dots, D$
and $j=1,2$. Use \eqref{eq1} and the assumption $\tau \circ
\psi_1 = \tau \circ \psi_2$ to find $k \ge k_0$ such that 
$\| \psi_j(f) - \psi_j^{(k)}(f)\| < \ep/3$ and
$$(\psi_j^{(k)})_*(h_{D,d-1}-h_{D,d}) \ge c, \qquad
\|(\psi_1^{(k)})_*(h_{D,d})-(\psi_2^{(k)})_*(h_{D,d})\|_\infty < c,
$$
for all $f \in F$, for $d=(1),2, \dots, D$, and for $j=1,2$. It then
follows from \cite[Corollary~5.6]{JiaSu:Z} that there is a unitary
element $u$ in $B_k$ such that
$\|\psi_2^{(k)}(f)-u^*\psi_1^{(k)}(f)u\| \le \ep/3$ for all $f \in F$;
whence $\|\psi_2(f)-u^*\psi_1(f)u\| \le \ep$ for all $f \in F$. This
proves that $\psi_1$ and $\psi_2$ are approximately unitarily
equivalent.
\end{proof}

\noindent
For any natural numbers $n$ and $m$ let $E(n,m)$ be the \Cs{} that
consists of all functions $f$ in 
$C([0,1], M_{n^\infty} \otimes M_{m^\infty})$ for which $f(0) \in
M_{n^\infty} \otimes \C$ and $f(1) 
\in \C \otimes M_{m^\infty}$. 

\begin{proposition} \label{prop:E(n,m)}
There is a unital embedding of $E(n,m)$ into $\cZ$ for every pair of
natural numbers $n,m$ that are relatively prime. 
\end{proposition} 

\begin{proof}
For each $k$ there is a unital embedding $\sigma_k
\colon M_{n^k} \otimes M_{m^k} \to M_{n^{k+1}} \otimes M_{m^{k+1}}$
which satisfies $$\sigma_k(M_{n^k} \otimes \C) \subseteq M_{n^{k+1}}
\otimes \C, \qquad \sigma_k(\C \otimes M_{m^k}) \subseteq \C \otimes
M_{m^{k+1}}.$$ 
Thus $f \mapsto \sigma_k
\circ f$ defines a \sh{} $\rho_k \colon I(n^k,n^k m^k,m^k) \to I(n^{k+1},
n^{k+1}m^{k+1}, m^{k+1})$, and  
$E(n,m)$ is the inductive limit of the sequence
$$\xymatrix{I(n,nm,m) \ar[r]^-{\rho_1} & I(n^2,n^2m^2,m^2)
  \ar[r]^-{\rho_2} & I(n^3,n^3m^3,m^3) \ar[r]^-{\rho_3} & \cdots
  \ar[r] & E(n,m).}
$$
Take standard unital embeddings $\psi_k \colon I(n^k,n^km^k,m^k) \to
\cZ$ (cf.\ Theorem~\ref{thm:I}~(i)). Then $\psi_k$ and  $\psi_{k+1} \circ
\rho_k$ are both standard unital embedding of $I(n^k,n^km^k,m^k)$
into $\cZ$, so they are approximately unitarily equivalent by 
Theorem~\ref{thm:I}~(ii). We
obtain  the desired embedding of $E(n,m)$ into $\cZ$ from this fact
combined with a one-sided approximate intertwining (after Elliott), see for
example \cite[Theorem~1.10.14]{Lin:amenable}.  
\end{proof}

\section{Almost unperforation} \label{sec:aup}

\noindent Consider an ordered abelian semigroup $(W, +, \le)$. An
element $x \in W$ is called \emph{positive} if $y + x \ge y$ for all
$y \in W$, and $W$ is said to be positive if all elements in $W$ are
positive. If $W$ has a zero-element $0$, then $W$ is positive if and
only if $0 \le x$ for all $x \in W$. An abelian semigroup equipped
with the \emph{algebraic order}: $x \le y$ iff $y = x+z$ for some $z
\in W$, is positive. 

\begin{definition} \label{def:aup}
A positive ordered abelian semigroup $W$ is said to be \emph{almost
  unperforated} if for all $x,y \in W$ and all $n,m \in \N$, with 
$nx \le my$ and $n > m$, one has $x \le y$.
\end{definition}
 
\noindent Let $W$ be a positive ordered abelian semigroup. 
Write $x \propto y$ if $x,y$ are elements in $W$ and $x \le ny$ for 
some natural number $n$ (i.e., $x$
belongs to the ideal in $W$ generated by the element
$y$). The element $y$ is said to be an \emph{order unit} for $W$ if $x
\propto y$ for all $x \in W$.  For each positive element $x$
in $W$ let $S(W,x)$ be the set of order preserving additive maps $f
\colon W \to [0,\infty]$ such that $f(x)=1$. Although we shall not use 
this fact, we mention that $S(W,x)$ is non-empty if and only if for all natural
numbers $n$ and $m$, with $nx \le mx$, one has $n \le m$. This follows
from \cite[Corollary~2.7]{BlaRor:extending} and the following observation
that also will be used in the proof of the proposition below. For any
element $x \in W$ the set $W_0=\{z \in W \mid z \propto x\}$ is an order 
ideal in $W$, and $x$ is an order unit for $W_0$. Moreover, any state
$f$ in $S(W_0,x)$ extends to a state $\overline{f}$ in $S(W,x)$
by setting  $\overline{f}(z) = \infty$ for $z \in W
\setminus W_0$. 

\begin{proposition} \label{prop:GH}
Let $W$ be a positive ordered abelian semigroup. Then $W$ is almost
unperforated if and only if the following condition holds: 
For all elements $x,y$ in $W$, with $x \propto y$ and $f(x) <
f(y)$ for all $f \in S(W,y)$, one has $x \le y$. 
\end{proposition}

\begin{proof} Following the argument above we can---if necessary by 
passing to an order ideal of $W$---assume that $y$ is an order unit
for $W$. The ``only if'' part now follows from
\cite[Pro\-po\-si\-tion~3.1]{Ror:UHFII}, which again uses Goodearl and
Handelman's extension result \cite[Lemma~4.1]{GooHan:extending}. 

To prove the ``if'' part, take elements $x,y \in W$ and $n \in \N$
such that $(n+1)x \le ny$. Then $x \propto y$ because $x \le (n+1)x \le ny$;
and $f(x) \le n(n+1)^{-1} < 1 = f(y)$ for
all $f \in S(W,y)$, whence $x \le y$. 
\end{proof} 

\begin{definition} \label{def:aup-G}
An ordered abelian group $(G,G^+)$ is said to be \emph{almost
  unperforated} if for all $g \in G$ and for all $n \in \N$, 
with $ng, (n+1)g \in G^+$, one has $g \in G^+$. 
\end{definition} 

\begin{lemma} \label{lm:G-G+}
Let $(G,G^+)$ be an ordered abelian group. Then $G$ is almost
unperforated if and only if the positive semigroup $G^+$ is almost
unperforated. 
\end{lemma}

\begin{proof}
Suppose that $G$ is almost unperforated and that $x,y \in G^+$ satisfy
$(n+1)x \le ny$ for some natural number $n$. Then $n(y-x) \ge x \ge 0$
and $(n+1)(y-x) \ge y \ge 0$, whence $y-x \ge 0$ and $y \ge
x$. Conversely, suppose that $G^+$ is almost unperforated and that
$ng, (n+1)g \in G^+$ for some $n \in \N$. Since $(n+1)ng = n(n+1)g$,
we get $ng \le (n+1)g$, which implies that $g = (n+1)g-ng$ belongs to $G^+$.
\end{proof}

\noindent A \emph{simple} ordered abelian group is almost unperforated if and
only if it is weakly unperforated. Indeed, if $n \in \N$ and $g \in G$
are such that $ng \in G^+ \setminus \{0\}$, then, by simplicity of $G$
there is a natural number $k$ such that $kng \ge g$. Thus
$(kn-1)g$ and $kng$ are positive, so $g$ is positive if $G$ is almost
unperforated (cf.\ Lemma~\ref{lm:G-G+}). Conversely, if $G$ is weakly
unperforated and $ng, (n+1)g \in G^+$, then $g \in G^+$ if $ng \ne
0$, and $g = (n+1)g \in G^+$ if $ng=0$. 

Elliott considered in \cite{Ell:torsion} a notion of what he called weak
unperforation of (non-simple) ordered abelian groups with torsion. (We
have refrained from using the term ``weak unperforation'' in 
Definition~\ref{def:aup-G} to avoid conflict with Elliott's definition.)
A torsion free
group is weakly unperforated in the sense of Elliott if and only if it
is unperforated: $ng \ge 0$ implies $g \ge 0$ for all group elements
$g$ and for all natural numbers $n$. The group $G = \Z^2$ with the
positive cone generated by the three elements $(1,0), (0,1), (2,-2)$
is torsion free with perforation: $(2,-2) \in G^+$ but $(1,-1) \notin
G^+$, so it is not weakly unperforated in the sense of Elliott. However, 
the group $(G,G^+)$ is almost unperforated, as the reader can verify. 

In the converse
direction, any weakly unperforated group is almost unperforated. Indeed, if 
$G$ is weakly unperforated and $g \in G$ and $n \in \N$ are such that 
$ng, (n+1)g$ are positive, then $g$ is positive modulo torsion, i.e., $g+t$
is positive for some $t \in G_{\mathrm{tor}}$. Let $k \in \N$ be the order of
$t$, find natural numbers $\ell_1, \ell_2$ such that 
$N=\ell_1n + \ell_2(n+1)$ is congruent with -1 modulo $k$. Then 
$Ng = Ng + (N+1)t$ is positive, whence $-t \le N(g+t)$, which by the 
hypothesis of weak unperforation implies that $g = (g+t) + (-t)$ is positive.

\section{Weak and almost unperforation of $\cZ$-absorbing \Cs s} \label{sec:1}

\noindent Cuntz associates in \cite{Cuntz:dimension} 
to each \Cs{} $A$ a positive ordered abelian
semigroup $W(A)$ as follows. Let $M_\infty(A)^+$ denote the (disjoint)
union $\bigcup_{n=1}^\infty M_n(A)^+$. For $a \in M_n(A)^+$ and $b \in
M_m(A)^+$ set $a \oplus b = \diag(a,b) \in M_{n+m}(A)^+$, 
and write $a \precsim b$ if there is a sequence $\{x_k\}$
in $M_{m,n}(A)$ such that $x_k^*bx_k \to a$. Write $a \sim b$ if $a
\precsim b$ and $b \precsim a$. Put $W(A) = M_\infty(A)^+/\! \sim$,
and let $\langle a \rangle  \in W(A)$ be the equivalence
class containing $a$ (so that $W(A) = 
\{ \langle a \rangle \mid a \in M_\infty(A)^+\}$). Then $W(A)$ is a positive
ordered abelian semigroup when equipped with the relations:
$$\langle a \rangle + \langle b \rangle = \langle a \oplus b \rangle,
\qquad \langle a \rangle \le \langle b \rangle \iff a \precsim b,
\qquad a,b \in M_\infty(A)^+.$$

Following the standard convention, for each positive element $a \in A$ and
for each $\ep \ge 0$, write $(a-\ep)_+$ for the positive
element in $A$ given by $h_\ep(a)$, where $h_\ep(t) = \max\{t-\ep,0\}$. 
We recall below some facts about
the comparison of two positive elements $a,b$ in a \Cs{} $A$ (see
\cite[Proposition~1.1]{Cuntz:dimension} and \cite[Section~2]{Ror:UHFII}):
\begin{itemize}
\item[(a)] $a \precsim b$
  if and only if $(a-\ep)_+ \precsim b$ for all $\ep>0$.
\item[(b)] $a \precsim b$ if and only if for each $\ep >0$ there
  exists $x \in A$ such that $x^*bx = (a-\ep)_+$.
\item[(c)] If $\|a-b\| < \ep$, then $(a-\ep)_+ \precsim b$.
\item[(d)] $\big((a-\ep_1)_+-\ep_2\big)_+ =
  \big(a-(\ep_1+\ep_2)\big)_+$.
\item[(e)] $a + b \precsim a \oplus b$; and if $a
  \perp b$, then $a + b \sim a \oplus b$.
\end{itemize}  
If $a$ belongs to the closed two-sided ideal, $\overline{AbA}$,
generated by $b$, then $(a-\ep)_+$ belongs to the algebraic two-sided
ideal, $AbA$, generated by $b$ for all $\ep>0$, in which case $(a-\ep)_+
= \sum_{i=1}^n 
x_i^*bx_i$ for some $n \in \N$ and some $x_i \in A$. This shows that
\begin{equation} \label{eq:ab}
a \in \overline{AbA} \iff \forall \ep > 0 \exists n \in \N: \langle
(a-\ep)_+ \rangle \le n \langle b \rangle.
\end{equation}

\begin{lemma} \label{lm:tensor}
Let $A$ and $B$ be two \Cs s, let $a,a' \in A$ and $b,b' \in B$ be
positive elements, and let $n,m$ be natural numbers.
\begin{enumerate}
\item If $n \langle a \rangle \le m \langle a' \rangle$ in $W(A)$, then
 $n \langle a \otimes b \rangle \le m \langle a' \otimes b \rangle$ in
 $W(A \otimes B)$.
\item If $n \langle b \rangle \le m \langle b' \rangle$ in $W(B)$, then
 $n \langle a \otimes b \rangle \le m \langle a \otimes b' \rangle$ in
 $W(A \otimes B)$.
\end{enumerate}
\end{lemma}

\begin{proof} (i). Assume that $n \langle a \rangle \le m \langle a'
  \rangle$ in $W(A)$. Then there is a sequence $x_k = \{x_k(i,j)\}$ in
  $M_{m,n}(A)$ such that $x_k^*(a' \otimes 1_m)x_k \to a \otimes 1_n$
  (or, equivalently, such that $\sum_{l=1}^m x_k(l,i)^*a'x_k(l,j) \to
  \delta_{ij} a$ for all $i,j = 1, \dots, n$).
  Let $\{e_k\}$ be a sequence of po\-si\-ti\-ve contractions in $B$ such
  that $e_kbe_k \to b$. Put $y_k(i,j) = x_k(i,j) \otimes e_k \in A
  \otimes B$, and put $y_k = \{y_k(i,j)\} \in M_{m,n}(A \otimes
  B)$. Then $y_k^*((a' \otimes b)
  \otimes 1_m)y_k \to (a \otimes b) \otimes 1_n$ (or, equivalently,
  $\sum_{l=1}^m y_k(l,i)^*(a' \otimes b)y_k(l,j) \to \delta_{ij} (a
\otimes b)$ for all $i,j = 1, \dots, n$). This shows that $n \langle a
\otimes b \rangle \le m \langle a' \otimes b \rangle$.

(ii) follows from (i) by symmetry.
\end{proof}

\begin{lemma} \label{lm:n-div}
For all natural numbers $n$ there exists a positive element $e_n$ in
$\cZ$ such that $n \langle e_n \rangle \le \langle 1_\cZ \rangle \le
(n+1) \langle e_n \rangle$. 
\end{lemma} 

\begin{proof} By Theorem~\ref{thm:I} (a fact which follows easily 
from from Jiang and Su's
  paper \cite{JiaSu:Z}) the \Cs{} $I=I(n,n(n+1),n+1)$ admits a unital
  embedding into $\cZ$, so it suffices to find a positive element
  $e_n$ in $I$ such that  $n \langle e_n \rangle \le \langle
  1_I \rangle \le (n+1) \langle e_n \rangle$ in $W(I)$. 

The idea of the proof is simple (but verifying the details requires some
effort): There are positive functions $f_1,
f_2, \dots, f_{n+1}$ in $I$ such that 
\begin{enumerate}
\item  $$f_i(0) = \begin{cases}  e_{ii}^{(n)} \otimes 1, & i=1, \dots, n,\\
0, & i=n+1, \end{cases} \qquad f_i(t) =  1 \otimes e_{ii}^{(n+1)}, \; \;
t \in [1/2,1],$$
(where
$\{e_{ij}^{(m)}\}_{i,j=1}^m$ denotes the canonical set of matrix units
for $M_m(\C)$), 
\item $f_1, f_2, \dots, f_n$ are pairwise orthogonal, 
\item $\sum_{i=1}^{n+1} f_i = 1$, and
\item $f_{n+1} \precsim f_1 \sim f_2 \sim \cdots \sim f_n$.  
\end{enumerate}
It will follow from (iii) and (e) that $\sum_{i=1}^{n+1} \langle f_i
\rangle \ge \langle 1 \rangle$; and (ii) and (e) imply that 
$\sum_{i=1}^{n} \langle f_i \rangle \le \langle 1 \rangle$. It
therefore follows from (iv) that $e_n=f_1$ has the desired property. 

We proceed to construct the functions $f_1, \dots, f_{n+1}$. Put
\begin{equation} \label{eq:W}
W = \sum_{i,j=1}^n  e_{ij}^{(n)} \otimes e_{ji}^{(n+1)} + 1 \otimes
e_{n+1,n+1}^{(n+1)}.
\end{equation}
Then $W$ is a self-adjoint unitary element in $M_n \otimes M_{n+1}$,
and 
\begin{equation} \label{eq:W1}
W(1 \otimes e_{ii}^{(n+1)})W^* =  e_{ii}^{(n)} \otimes (1 -
e_{n+1,n+1}^{(n+1)}) \le e_{ii}^{(n)} \otimes 1,
\end{equation}
for $1 \le i \le n$. Choose a continuous path of unitaries $t \mapsto
V_t$ in $M_n \otimes M_{n+1}$, $t \in [0,1]$, such that $V_0 = 1$ and
$V_t = W$ for $t \in [1/2,1]$. Put $W_t = V_tW$. Choose a continuous
path $t \mapsto \gamma_t \in [0,1]$ such that $\gamma_0
= 0$ and $\gamma_t = 1$ for $t \in [1/2,1]$. Define
$f_i \colon [0,1] \to M_n \otimes M_{n+1}$ by 
\begin{eqnarray*} 
f_i(t) & = &  \gamma_t W_t(1 \otimes e_{ii}^{(n+1)})W_t^* +
(1-\gamma_t)V_t(e_{ii}^{(n)} \otimes 1)V_t^*, \qquad 
i=1, \dots, n,\\  
f_{n+1}(t) & = & \gamma_t W_t ( 1 \otimes e_{n+1,n+1}^{(n+1)})W_t^*,
\end{eqnarray*}
where $t \in [0,1]$. It is easy to check that (i) and (iii) above hold.
From (i) we see that all $f_i$ belong to $I$. Use
\eqref{eq:W1} to see that $f_i(t) \le V_t(e_{ii}^{(n)} \otimes
1)V_t^*$ for all $t$ and all $i=1, \dots, n$, and use again this to see that 
(ii) holds. 

We proceed to show that (iv) holds. Let $S \in M_n$ and $T \in M_{n+1}$ 
be the permutation unitaries for which 
$S^{i-1}e_{11}^{(n)}S^{-(i-1)} = e_{ii}^{(n)}$ and
$T^{i-1}e_{11}^{(n+1)}T^{-(i-1)} = e_{ii}^{(n+1)}$ for 
$i = 1, \dots, n, (n+1)$. Put
$$R_i(t) = V_t (S^{i-1} \otimes 1)V_t^*, \qquad t \in [0,1/2], \; i
=1, \dots, n.$$
Brief calculations show that $R_i(0) = S^{i-1} \otimes 1 \in M_n \otimes
\C$, $f_i(t) = R_i(t)f_1(t)R_i(t)^*$ for $t
\in [0,1/2]$,  and 
\begin{eqnarray*}
R_i(1/2)(1 \otimes e_{11}^{(n+1)}) R_i(1/2)^* & = &
R_i(1/2)f_1(1/2)R_i(1/2)^* \, = \, f_i(1/2) \, = \,
1 \otimes e_{ii}^{(n+1)} 
\\ & =& (1 \otimes T^{(i-1)})(1 \otimes e_{11}^{(n+1)})(1 \otimes
T^{-(i-1)}),
\end{eqnarray*}
for $i=1, \dots, n$. The unitary group of the relative commutant
$M_n \otimes M_{n+1} \cap \{1 \otimes e_{11}^{(n+1)}\}'$ is connected, 
so we can 
extend the paths $t \mapsto R_i(t)$, $t \in [0,1/2]$, to continuous 
path $t \mapsto R_i(t)$, $t \in [0,1]$, such 
that $R_i(1) = 1 \otimes T^{-(i-1)} \in \C \otimes M_{n+1}$ and
$R_i(t)f_1(t)R_i(t)^* =  1 \otimes e_{ii}^{(n+1)} = f_i(t)$ for $t
\in [1/2,1]$ and for $i=2, \dots, n$. Thus $R_i$ is a unitary element in
$I$ and $R_i f_1 R_i^* = f_i$ for  $i=2, \dots, n$. This proves that $f_1 \sim
f_2 \sim \cdots \sim f_n$. 

We must also show that $f_{n+1} \precsim f_1$. Let $g_i
\in I$ be given by $g_i(t) = \gamma_t W_t(1 \otimes
e_{ii}^{(n+1)})W_t^*$ for $i = 1, \dots, n+1$ (so that $f_{n+1} = g_{n+1}$). 
A calculation then shows that $R_ig_1R_i^* = g_i$ so that $g_1, \dots, g_n$ 
are unitarily equivalent. 
By symmetry, $f_{n+1} = g_{n+1}$ is unitarily equivalent
to $g_1$, and as $g_1 \le f_1$, we conclude that $f_{n+1} \precsim
f_1$. 
\end{proof}

\begin{lemma} \label{lm:A_k}
Let $A$ be any \Cs, and let $a,a'$ be positive
elements in $A$ for which $(n+1) \langle a \rangle \le n
\langle a' \rangle$ in $W(A)$ for some natural number $n$. Then
$\langle a \otimes 1_\cZ \rangle \le \langle a' \otimes 1_\cZ \rangle$
in $W(A \otimes \cZ)$.
\end{lemma}

\begin{proof} Take $e_n$ in $\cZ$
  as in Lemma~\ref{lm:n-div}. Then, by Lemma~\ref{lm:tensor}, 
$$\langle a \otimes 1_\cZ \rangle \le (n+1) \langle a \otimes e_n
\rangle \le n \langle a' \otimes e_n \rangle \le \langle a' \otimes
1_\cZ \rangle$$
in $W(A \otimes \cZ)$.
\end{proof}

\begin{lemma} \label{lm:decomposition}
Let $A$ be a $\cZ$-absorbing \Cs{}. Then there is a sequence of
isomorphisms $\sigma_n \colon A \otimes \cZ \to A$ such that 
$$\lim_{n \to \infty} \|\sigma_n(a \otimes 1) - a\| = 0, \qquad a \in A.$$
\end{lemma}

\begin{proof} It is shown in \cite{JiaSu:Z} that $\cZ$ is isomorphic
  to $\bigotimes_{k=1}^\infty \cZ$. We may therefore identify $A$ with
  $A \otimes (\bigotimes_{k=1}^\infty \cZ)$. With this identification
we define 
$\sigma_n \colon A \otimes (\bigotimes_{k=1}^\infty \cZ) \otimes \cZ  \to A
\otimes (\bigotimes_{k=1}^\infty \cZ)$
to be the isomorphism that fixes $A$ and the first $n$ copies of
$\cZ$, which sends the last copy of $\cZ$ to the copy of $\cZ$ at
position ``$n+1$'', and which shifts the remaining copies of
$\cZ$ one place to the right. 
\end{proof}

\begin{theorem} \label{thm:wup}
Let $A$ be a $\cZ$-absorbing \Cs{}. Then $W(A)$ is almost unperforated.
\end{theorem}

\begin{proof}
Let $a,a'$ be positive elements in $M_\infty(A)$ for which
$(n+1)\langle a \rangle \le n \langle a' \rangle$. Upon replacing $A$
by a matrix algebra over $A$ (which still is $\cZ$-absorbing) we may
assume that $a$ and $a'$ both belong to $A$. Let $\ep >0$. It follows
from Lemma~\ref{lm:A_k} that $a \otimes 1 \precsim a' \otimes 1$ in $A
\otimes \cZ$, so there exists $x \in A \otimes \cZ$ with $\|x^*(a'
\otimes 1)x - a \otimes 1\| < \ep$. 
Let $\sigma_k \colon A \otimes \cZ \to A$
be as in Lemma~\ref{lm:decomposition}, and put $x_k =
\sigma_k(x)$. Then $\|x_k^*\sigma_k(a' \otimes 1)x_k - \sigma_k(a \otimes 1)\| 
< \ep$, whence $\|x_k^*a'x_k - a\| < \ep$ if $k$ is chosen large enough. 
This shows that 
$\langle a \rangle \le \langle a' \rangle$ in $W(A)$.
\end{proof}

\noindent A \Cs{} $A$ where $W(A)$ is almost unperforated has nice
comparability properties, as we shall proceed to illustrate in the
remaining part of this section, and in the later sections of this
paper. 

We recall a few facts about \emph{dimension function},
introduced by Cuntz in \cite{Cuntz:dimension}. A dimension function 
on a \Cs{} $A$ is an additive order preserving function $d \colon W(A) \to
[0,\infty]$. (We can also regard  
$d$ as a function $M_\infty(A)^+ \to [0,\infty]$ that respects the
rules $d(a \oplus b) = d(a) + d(b)$ and $a \precsim b \Rightarrow d(a)
\le d(b)$ for all $a,b \in M_\infty(A)^+$.) A dimension function $d$
is said to be \emph{lower semi-continuous} if $d = \overline{d}$, where
\begin{equation} \label{eq:d-bar}
\overline{d}(a) \eqdef \lim_{\ep \to 0+} d((a-\ep)_+), \qquad a \in
M_\infty(A)^+.
\end{equation}
Moreover, $\overline{d}$ is a lower semi-continuous dimension function on 
$A$ for each dimension function $d$, cf.\ \cite[Proposition~4.1]{Ror:UHFII}. 
Note that 
$d((a-\ep)_+) \le \overline{d}(a) \le d(a)$ for every dimension
function $d$ and for every $\ep >0$, and that $\overline{d}(p) = d(p)$
for every projection $p$. 

By an \emph{extended trace} on a \Cs{} $A$ we shall mean a function
$\tau \colon A^+ \to [0,\infty]$ which is additive, homogeneous, and
has the trace property: $\tau(x^*x)=\tau(xx^*)$ for all $x \in A$. 
 If $\tau$ is an extended trace on $A$, then 
\begin{equation} \label{eq:d-tau}
d_\tau(a) \eqdef \lim_{\ep \to 0+} \tau(f_\ep(a)) \, \, (= \lim_{n \to
  \infty} \tau(a^{1/n})), \quad a \in M_\infty(A)^+,
\end{equation}
where $f_\ep \colon \R^+ \to \R^+$ is given by $f_\ep(t) =
\min\{\ep^{-1}t,1\}$, 
defines a lower semi-continuous dimension function on $A$.
If $A$ is exact, then every lower semi-continuous dimension function
on $A$ is of the form $d_\tau$ for some extended trace $\tau$ on $A$.
(This follows from Blackadar and Handelman,
\cite[Theorem~II.2.2]{BlaHan:quasitrace}, who show
that one can lift $d$ to a quasitrace, and from Haagerup, \cite{Haa:quasi},
and Kirchberg, \cite{Kir:quasitraces}, who show that
quasitraces are traces on exact \Cs s).  

With the characterization of lower semi-continuous dimension functions
above and with Theorem~\ref{thm:wup} at hand one can apply the proof of
\cite[Theorem~5.2]{Ror:UHFII} to obtain:

\begin{corollary} \label{cor:compare-simple}
Let $A$ be a \Cs{} for which $W(A)$ is almost unperforated (in
particular, $A$ could be a $\cZ$-absorbing \Cs), and suppose in
addition that $A$ is exact, simple and unital. Let $a,b$
be positive elements in $A$. If $d_\tau(a) < d_\tau(b)$ for every
tracial state $\tau$ on $A$, then $a \precsim b$.
\end{corollary}  

\noindent We also have the following ``non-simple'' version of the
result above. 

\begin{corollary} \label{cor:compare}
Let $A$ be a \Cs{} for which $W(A)$ is almost unperforated (in
particular, $A$ could be a $\cZ$-absorbing \Cs). Let $a,b$
be positive elements in $A$. Suppose that $a$ belongs to
$\overline{AbA}$ and that $d(a) < d(b)$ for every
dimension function $d$ on $A$ with $
d(b) =1$. Then $a \precsim b$.  
\end{corollary}

\begin{proof}
It follows from \eqref{eq:ab} that  $\langle (a-\ep)_+ \rangle \propto
\langle b \rangle$ in $W(A)$ for each $\ep >0$; and by assumption,
$d(\langle (a-\ep)_+ \rangle) \le d(\langle a \rangle) 
< d(\langle b \rangle)$ for every $d \in
S(W(A),\langle b \rangle)$. Thus 
$\langle (a-\ep)_+ \rangle \le \langle b \rangle$ by
Proposition~\ref{prop:GH}, and this proves the corollary as $\ep>0$ was
arbitrary.  
\end{proof} 

\noindent Gong, Jiang and Su proved in \cite{GongJiangSu:Z} that the
$K_0$-group of a simple unital $\cZ$-absorbing \Cs{} is weakly
unperforated. At the level of semigroups, we can extend this result to
the non-simple case, as explained below.

Let $V(A)$ denote the semigroup of Murray-von Neumann
equivalence classes of projections in matrix algebras over $A$
equipped with the algebraic order: $x \le y$ if there exists $z$ such
that $y = x+z$. The relation ``$\precsim$'', defined in the beginning 
of this section, agrees with the usual comparison relation when applied
to projections $p$ and $q$, i.e.,  $p \precsim q$ if and only if $p$ is 
equivalent to a subprojection of $q$. The corollary below is thus an 
immediate consequence of Theorem~\ref{thm:wup}.  

\begin{corollary} \label{cor:V(A)}
The semigroup $V(A)$ is almost unperforated for every $\cZ$-absorbing
\Cs{} $A$. 
\end{corollary}

\noindent It follows from Lemma~\ref{lm:G-G+} that if
$A$ is a stably finite \Cs{} with an approximate unit consisting of
projections, then $K_0(A)$ is almost unperforated if and 
only if $K_0(A)^+$ is almost unperforated. It seems plausible that
$K_0(A)^+$ is almost unperforated whenever $V(A)$ is almost
unperforated; and this is trivially the case when $V(A)$ has the cancellation
property. This implication also holds when
$V(A)$ is simple. Indeed, let $\gamma \colon V(A) \to K_0(A)$ be the
Grothendieck map, so that $K_0(A)^+ = \gamma(V(A))$. Take $x,y \in
V(A)$ and $n \in \N$ such that $(n+1)\gamma(x) \le n \gamma(y)$. Then
$(n+1)x +u \le ny +u$ for some $u \in V(A)$. Repeated use of this
inequality yields $N(n+1)x + u \le Nny + u$ for all natural numbers
$N$. That $V(A)$ is simple means that every non-zero element, and
hence $y$, is an order unit for $V(A)$, so there is a natural number 
$k$ with $u \le ky$. Now, $N(n+1)x \le N(n+1)x+u \le Nny +u \le
(Nn+k)y$, which for $N \ge k+1$ yields $x \le y$ and hence $\gamma(x)
\le \gamma(y)$.  

We thus have the following result, that
slightly extends \cite[Theorem~1]{GongJiangSu:Z}. 

\begin{corollary} \label{cor:K0-wup}
Let $A$ be a stably finite $\cZ$-absorbing \Cs{} with an approximate
unit consisting of projections. If $V(A)$ has the cancellation
property or if $V(A)$ is simple, then $K_0(A)$ is almost unperforated.
\end{corollary}

\begin{corollary} \label{cor:compare-proj}
Let $A$ be an exact \Cs{} for which $W(A)$ is almost unperforated (in
particular, $A$ could be an exact $\cZ$-absorbing \Cs). Let $p,q$
be projections in $A$ such that $p$ belongs to
$\overline{AqA}$. Suppose that $\tau(p) < \tau(q)$ for every
extended trace $\tau$ on $A$ with $\tau(q)=1$. Then $p \precsim q$.  
\end{corollary}

\begin{proof}
We show that $d(p) < d(q)$ for every dimension function $d$ on $A$
with $d(q)=1$,
and the result will then follow from Corollary~\ref{cor:compare}. Let
$\overline{d}$ be the lower semi-continuous dimension function  
associated with $d$ in \eqref{eq:d-bar}.
As remarked above, by Haagerup's theorem on quasitraces,  
$\overline{d} = d_\tau$ for some extended trace $\tau$, cf.\ \eqref{eq:d-tau}.
Because $d$,
$\overline{d}$ and $\tau$ agree on projections, we have
$\tau(q)=d(q)=1$ and  $d(p) = \tau(p) < \tau(q) = d(q)$ as
desired.
\end{proof}

\section{Applications to purely infinite \Cs s}

\noindent In this short section we derive two results that say 
when a $\cZ$-absorbing \Cs{} is purely infinite and
$\cO_\infty$-absorbing. Similar results were obtained in
\cite{KirRor:pi2} for approximately divisible \Cs s. An exact \Cs{} 
is said to be \emph{traceless} when it admits no
extended trace (see Section~3) that takes values other than $0$ and $\infty$. 

\begin{corollary} \label{cor:pi}
Let $A$ be an exact \Cs{} for which $W(A)$ is almost unperforated (in
particular, $A$ could be an exact $\cZ$-absorbing \Cs). Then $A$ is
purely infinite 
if and only if $A$ is traceless. 
\end{corollary}

\begin{proof} Note first that $A$, being traceless, can have no
  abelian quotients. Each lower semi-continuous
  dimension function on $A$ arises from an extended trace on $A$ (as
  remarked above Corollary~\ref{cor:compare-simple}), and must therefore
  take values in $\{0,\infty\}$, again because $A$ is traceless. 

  Take positive elements $a,b$ in $A$ such that 
 $a$ belongs to $\overline{AbA}$. We must show that $a \precsim b$
 (cf.\ \cite{KirRor:pi}), and it suffices to show that $(a-\ep)_+
 \precsim b$ for all $\ep>0$. Take $\ep>0$. Let $d$ be a
 dimension function on $A$ such that $d(b)=1$ (if such a dimension
 function exists), and let $\overline{d}$ be its associated lower
 semi-continuous dimension function, cf.\ \eqref{eq:d-bar}. 
 Then $\overline{d}(b) = 0$ (as remarked
 above). Use \eqref{eq:ab} to see that $\overline{d}((a-\ep/2)_+) =
 0$, and hence that $d((a-\ep)_+)=0$. Thus $(a-\ep)_+ \precsim b$ by
 Corollary~\ref{cor:compare}.
\end{proof}

\begin{theorem} \label{thm:spi}
Let $A$ be a nuclear separable $\cZ$-absorbing \Cs. Then $A$ absorbs
$\cO_\infty$ (i.e., $A \cong A \otimes \cO_\infty$) if and only if $A$
is traceless.
\end{theorem}

\begin{proof}
If $A$ absorbs $\cO_\infty$, then $A$ is traceless (see
\cite[Theorem~9.1]{KirRor:pi2}). Suppose conversely that $A$ is traceless.
We then know from Corollary~\ref{cor:pi} that $A$
is purely infinite. It follows from \cite[Theorem~9.1]{KirRor:pi2} (in
the unital or the stable case) and from \cite[Corollary~8.1]{Kir:spi} (in 
the general case) that $A$ absorbs $\cO_\infty$ if $A$ is \emph{strongly purely}
infinite, cf.\ \cite[Definition~5.1]{KirRor:pi2}. 

We need therefore only show that any $\cZ$-absorbing purely infinite \Cs{}
is strongly purely infinite. Let  
$$\begin{pmatrix} a & x \\ x^* & b \end{pmatrix} \in M_2(A)^+,$$
and let $\ep >0$ be given. Take pairwise orthogonal non-zero positive
contractions $h_1,h_2$ in the simple \Cs{} $\cZ$. Then
$$a \otimes 1 \in \overline{(A \otimes \cZ)(a \otimes h_1)(A \otimes
  \cZ)}, \qquad 
b \otimes 1 \in \overline{(A \otimes \cZ)(b \otimes h_2)(A \otimes
  \cZ)}.$$
Since $A \otimes \cZ \cong A$ is purely infinite there are elements
$c_1,c_2$ in $A \otimes \cZ$ such that 
$$\|c_1^*(a \otimes h_1)c_1 - a \otimes 1\| < \ep, \qquad
\|c_2^*(b \otimes h_2)c_2 - b \otimes 1\| < \ep.$$ Let $\sigma_n \colon
A \otimes \cZ \to A$ be as in Lemma~\ref{lm:decomposition}, and put
$d_{1,n} = \sigma_n((1 \otimes h_1^{1/2})c_1)$ and 
$d_{2,n}=\sigma_n((1 \otimes h_2^{1/2})c_2)$. Then 
\begin{eqnarray*}
& \|d_{1,n}^*\sigma_n(a \otimes 1)d_{1,n} - \sigma_n(a \otimes
1)\| < \ep, \qquad 
\|d_{2,n}^*\sigma_n(b \otimes 1)d_{2,n} - \sigma_n(b \otimes 1)\| < \ep,&\\
& d_{2,n}^*\sigma_n(x \otimes 1)d_{1,n} = 0.&
\end{eqnarray*}
The norm of $d_{j,n}$ does not dependent on $n$. Thus, if we take 
$d_1 = d_{1,n}$ and $d_2 = d_{2,n}$ for some large enough $n$, then we obtain
the desired estimates: $\|d_1^*ad_1-a\| <
\ep$, $\|d_2^*bd_2-b\|< \ep$, and $\|d_2^*xd_1\| < \ep$.
\end{proof} 

\section{The stable rank of $\cZ$-absorbing \Cs s}

\noindent We shall in this section show that simple, finite
$\cZ$-absorbing \Cs s have stable rank one. 

\begin{definition} A unital \Cs{} $A$ is said to be \emph{strongly
    $K_1$-surjective} if the canonical mapping $\cU(A_0 + \C 1_A) 
\to K_1(A)$ is surjective for every full hereditary sub-\Cs{} $A_0$ of $A$.
If the canonical mapping 
  $\cU(A_0 + \C 1_A)/\cU_0(A_0 + \C 1_A) \to K_1(A)$ is injective for 
every full hereditary sub-\Cs{} $A_0$, then we say that $A$ is
\emph{strongly $K_1$-injective}.
\end{definition}

\noindent Note that we do not assume simplicity in the two next lemmas. 

\begin{lemma} \label{lm:sp}
Every full hereditary sub-\Cs{} in a unital approximately divisible
\Cs{} contains a full projection.
\end{lemma}

\begin{proof} Let $B$ be a full hereditary sub-\Cs{} of a unital
  approximately divisible \Cs{} $A$. Take a full positive element $b$
  in $B$. Then $n\langle b \rangle \ge \langle 1_A \rangle$ in $W(A)$
  for some natural number $n$ (by \eqref{eq:ab}). Since $A$ is
  approximately divisible, 
there is a unital embedding of $M_{n+1} \oplus M_{n+2}$ into $A$, and,
as shown in \cite{TomsWin:Z}, $A$ is $\cZ$-absorbing. (We shall only apply this
lemma in the case where $A$ is the tensor product of a unital \Cs{}
with a UHF-algebra, and in this case we can conclude that $A$ is
$\cZ$-absorbing by the result in \cite{JiaSu:Z} that any
non-elementary simple AF-algebra, and in particular, every
UHF-algebra, is $\cZ$-absorbing.) 

Let  $e$ and $f$ be one-dimensional projections in $M_{n+1}$ and 
$M_{n+2}$, respectively, and let $p \in A$ be the image of $(e,f)$ 
under the inclusion mapping $M_{n+1} \oplus M_{n+2} \to A$. Then $p$ is a 
full projection that satisfies $(n+1)\langle p \rangle \le \langle 1
\rangle \le n \langle b \rangle$. Hence, by Theorem~\ref{thm:wup}, 
$\langle p \rangle \le \langle b \rangle$, i.e., $p \precsim b$. It
follows that $p = x^*bx$ for some $x \in A$. Put $v = b^{1/2}x$. Then
$p=v^*v$ and $p \sim vv^* = b^{1/2} xx^*b^{1/2} \in B$, so $vv^*$ is a full
projection in $B$.
\end{proof}

\begin{lemma} \label{lm:K1} Every unital approximately divisible \Cs{} is 
strongly $K_1$-surjective.
\end{lemma}

\begin{proof}
Let $B$ be a full hereditary sub-\Cs{} of $A$. We must show that the
canonical map $\cU(B + \C 1_A) \to K_1(A)$ is surjective. Use
Lemma~\ref{lm:sp} to find a full projection $p$ in $B$. It
suffices to show that the canonical map $\cU(pAp + \C(1_A-p)) \to
K_1(A)$ is surjective. Take an element $g$ in $K_1(A)$, and represent
$g$ as the class of a unitary element $u$ in $M_n(A)$ for some large
enough natural number $n$. Upon replacing $M_n(A)$ by $A$ we can
assume that $n=1$. 

Let $\cP$ be the set of projections $q \in A$ such that there exists a
unitary element $v \in qAq$ for which $g = [v + (1_A-q)]_1$ in
$K_1(A)$. We must show that $\cP$ contains all full projections in
$A$. Note first that if $q_1, q_2$ are projections in $A$ with $q_1
\precsim q_2$ and $q_1 \in \cP$, then $q_2 \in \cP$. Indeed, if $v_1
\in q_1Aq_1$ is unitary with $g = [v_1+(1_A -q_1)]_1$ and if $s^*s
= q_1$, $ss^* \le q_2$, then $v_2$ given by $sv_1s^* +(q_2-ss^*)$ is a
unitary element in $q_2Aq_2$, and $[v_1 + (1_A-q_1)]_1 = [v_2 +
(1_A-q_2)]_1$. 

Let $p \in A$ be a full projection. Then $(n-1) \langle p \rangle \ge
\langle 1_A \rangle$ in $W(A)$ for some large enough natural number
$n$, cf.\ \eqref{eq:ab}.
By approximate divisibility of $A$ there is a unitary element
$u_0 \in A$ with $\|u-u_0\| < 2$, and a unital embedding
$M_n \oplus M_{n+1} \to A \cap \{u_0\}'$; in other words, there are
matrix units $\{e_{ij}\}_{i,j=1}^n$ and $\{f_{ij}\}_{i,j=1}^{n+1}$ in
$A \cap \{u_0\}'$ such that $\sum_i e_{ii} + \sum_i f_{ii} =
1_A$. Note that $g = [u_0]_1$. Put $q = e_{11} + f_{11}$ and put
$v = u_0^ne_{11} + u_0^{n+1}f_{11}$, so that $v$ is a unitary element in
$qAq$. It follows from the Whitehead lemma that $u_0$ is homotopic to
$v+(1_A-q)$, whence $q$ belongs to $\cP$. As $n \langle q \rangle \le
\langle 1_A \rangle \le (n-1) \langle p \rangle$, it follows from
Theorem~\ref{thm:wup} that $q \precsim p$, whence $p \in \cP$ by the result
in the second paragraph of the proof. 
\end{proof}

\noindent Rieffel proved in \cite{Rfl:sr} that if $A$ is a unital
\Cs{} of stable rank one, then the canonical map $\cU(A)/\cU_0(A) \to
K_1(A)$ is an isomorphism, and hence injective. Rieffel also showed
for any such \Cs{} $A$ and any hereditary sub-\Cs{} $B$ of $A$
(full or not) that the stable rank of $B + \C 1_A$ is one. 

Take now a full hereditary sub-\Cs{} $B$ of $A$, where $A$ is unital
and of stable rank one. Then 
$$\cU(B + \C 1_A)/\cU_0(B + \C 1_A) \to K_1(B) \to K_1(A)$$ 
is an isomorphism (the second map is an isomorphism by Brown's
theorem, which guarantees that $A \otimes \cK \cong B \otimes
\cK$). 

This shows that any unital \Cs{} of stable rank one is strongly 
$K_1$-injective.

For each element $x$ in a \Cs{} $A$ we can write $x = v|x|$, where $v$
is a partial isometry in $A^{**}$. The element $x_\ep \eqdef
v(|x|-\ep)_+$ belongs to $A$ for every $\ep \ge 0$, and $\|x-x_\ep\|
\le \ep$. If $x$ is positive, then $x_\ep = (x-\ep)_+$. 

\begin{lemma} \label{lm:g_ep}
Let $A$ be a unital \Cs{}, let $a$ be a positive element
  in $A$, let $0 < \ep' < \ep$ be given, and set $A^0 =
  \overline{g_{\ep}(a)Ag_{\ep}(a)}$, where $g_{\ep} \colon \R^+ \to \R^+$
  is given by $g_{\ep}(t) = \max\{1-t/\ep,0\}$. Then $wa_\ep =
  a_\ep$ for every $w \in A^0 + 1_A$; and if $w$ is a unitary element
  in $A$ that satisfies $wa_{\ep'} = a_{\ep'}$, then $w$ belongs to
  $A^0 + 1_A$.
\end{lemma}

\begin{proof} The first claim follows from the fact that $xa_\ep = 0$
  for every $x \in A^0$. Suppose that $w$ is a unitary element in $A$
  with $wa_{\ep'} = a_{\ep'}$. Then $a_{\ep'}w = a_{\ep'}$, so $w-1_A$
  is orthogonal to $a_{\ep'}$. But the orthogonal complement of
  $a_{\ep'}$ is contained in $A^0$.
\end{proof}

\begin{proposition} \label{prop:pull-back}
Given a pull-back diagram
\begin{equation} \label{eq:pull-back}
\begin{split}
\xymatrix{& A \ar[dl]_-{\varphi_1} \ar[dr]^-{\varphi_2} \ar[dd]_-\pi& 
\\ A_1 \ar[dr]_-{\psi_1} & & A_2 \ar[dl]^-{\psi_2} \\ & B & }
\end{split}
\end{equation} 
with surjective \sh s $\psi_1$ and $\psi_2$.
Suppose that $A, A_1, A_2$ and $B$ are
unital \Cs s and that $a \in A$ are such that
\begin{enumerate}
\item $A_1$ and $A_2$ are strongly $K_1$-surjective,
\item $B$ is strongly $K_1$-injective,
\item $\image(K_1(\psi_1)) + \image(K_1(\psi_2)) = K_1(B)$,
\item $a^*a$ is non-invertible in every non-zero quotient of $A$. 
\end{enumerate}
Then $a$ belongs to the closure of
$\GL(A)$ if and only if $\varphi_j(a)$ belongs to the closure of
$\GL(A_j)$ for $j=1,2$.
\end{proposition}

\noindent The pull-back diagram \eqref{eq:pull-back} can, given $\psi_j \colon
A_j \to B$, $j=1,2$, be realized with
$A = \{(a_1,a_2) \in A_1 \oplus A_2 \mid \psi_1(a_1)=\psi_2(a_2)\}$ and
with $\varphi_j(a_1,a_2)=a_j$.

\begin{proof} The ``only if'' part is trivial. Assume now that
  $\varphi_j(a)$ belongs to the closure of $\GL(A_j)$ for $j=1,2$. Let
  $\ep>0$ be given. It
  then follows from \cite[Theorem~2.2]{Ror:unitary} that there are
  unitary elements $u_j$ in $A_j$ such that
  $\varphi_j(a_{\ep/2}) = u_j|\varphi_j(a_{\ep/2})|$ for $j=1,2$.
  We show below that there
  are unitary elements $v_j$ in $A_j$, $j=1,2$, such that $\varphi_j(a_\ep) =
  v_j|\varphi_j(a_\ep)|$, $j=1,2$, and $\psi_1(v_1)=\psi_2(v_2)$. It
  follows that $v=(v_1,v_2)$ is a unitary element in $A$ and that
  $a_\ep = v|a_\ep|$. This shows that $a$ belongs to the closure
  of the invertibles in $A$ (because $\ep>0$ was arbitrary). 

Let $g_{\ep} \colon \R^+ \to
\R^+$ be as in Lemma~\ref{lm:g_ep}, and put $A^0 = \overline{g_{\ep}(|a|)A
g_{\ep}(|a|)}$. Put
$$A_j^0 = \varphi_j(A^0) =
\overline{g_\ep(|\varphi_j(a)|) A_j g_\ep(|\varphi_j(a)|)}, \qquad
B^0 = \pi(A^0) = \overline{g_\ep(|\pi(a)|) B g_\ep(|\pi(a)|)}.$$
Assumption (iv) implies that $A^0$ is full in $A$. It follows that
the hereditary subalgebras $A_1^0, A_2^0, B^0$ are full in $A_1, A_2$
and $B$, respectively.

It follows from the identity
$$\pi(a_{\ep/2}) = \psi_1(u_1)|\pi(a_{\ep/2})|= 
\psi_2(u_2)|\pi(a_{\ep/2})|,$$ 
that $\psi_2(u_2)^*\psi_1(u_1)|\pi(a_{\ep/2})| = 
|\pi(a_{\ep/2})|$, and so $z \eqdef \psi_2(u_2)^*\psi_1(u_1)$ belongs to 
$B^0 + 1_B$ (cf.\ Lemma~\ref{lm:g_ep}). We show below that
$z=\psi_2(w_2)\psi_1(w_1^*)$ for some unitaries $w_j$ in 
$A_j^0 + \C 1_{A_j}$, $j=1,2$. 

Use conditions (i) and (iii) to find unitaries $y_j \in A_j^0 + \C
1_{A_j}$ such that $[\psi_2(y_2)\psi_1(y_1)^*]_1 = [z]_1$ in
$K_1(B)$. By condition (ii), the unitary element $(z_0 =) \,
z\psi_1(y_1)\psi_2(y_2^*)$ is homotopic to 1 in the unitary group of
$B^0 + \C 1_B$. Hence $z_0 = \psi_2(y_0)$ for some unitary $y_0$ in
$A_2^0 + \C 1_{A_2}$. Now, $w_1=y_1$ and $w_2 = y_0y_2$ are as desired.

Upon
replacing $w_1$ and $w_2$ by $\lambda w_1$ and $\lambda w_2$
for a suitable $\lambda \in \C$ with $|\lambda|=1$, we can assume that
$w_j \in A_j^0 + 1_{A_j}$. Then, by Lemma~\ref{lm:g_ep}, 
$w_j|\varphi_j(a_\ep)| = |\varphi_j(a_\ep)|$, $j=1,2$. It follows that
$(v_j = ) \; u_jw_j$ is a unitary in $A_j$, that $v_j|\varphi_j(a_\ep)| =
u_j|\varphi_j(a_\ep)| = \varphi_j(a_\ep)$ for $j=1,2$, and
$\psi_1(v_1) = \psi_1(u_1)\psi_1(w_1) = \psi_2(u_2)\psi_2(w_2) =
\psi_2(v_2)$, as desired.
\end{proof}


\begin{lemma} \label{lm:cts_field}
Let $A$ be a simple, unital, finite \Cs{}. Then $a \otimes 1$ belongs
to the closure of the invertibles in $A \otimes \cZ$ for every $a \in
A$.  
\end{lemma}

\begin{proof} Let $E_{2,3}$  be the \Cs{} which in
  Proposition~\ref{prop:E(n,m)} is shown to have a unital embedding
  into $\cZ$. It suffices to show that $a \otimes 1$ belongs to the
  invertibles in $A \otimes E_{2,3}$. 
If $a^*a \otimes 1$ is invertible in some non-zero quotient 
of $A \otimes E_{2,3}$, then $a^*a$ is invertible in $A$ by simplicity of $A$,
which again implies that $a$ is invertible, because $A$ is finite. The
claim of the lemma is trivial in this case. Suppose now that there is no
non-zero quotient of $A \otimes E_{2,3}$ in which $a^*a \otimes 1$ is 
invertible.

Identify $A \otimes E_{2,3}$ with 
$$\{f \in C([0,1],A \otimes B) \mid f(0) \in A \otimes B_1, 
\, \, f(1) \in A \otimes B_2 \},$$
where $B, B_1, B_2$ are UHF algebras of type $6^\infty$,
$2^\infty$, and $3^\infty$, respectively, with $B_j \subseteq B$. 
We have a pull-back diagram
\begin{equation} \label{eq:pull-back_2}
\begin{split}
\xymatrix{& A \otimes E_{2,3}
\ar[dl]_-{\varphi_1} \ar[dr]^-{\varphi_2} \ar[dd]_-{\ev_{\tfrac{1}{2}}}& 
\\ D_1 \ar[dr]_-{\ev_{\tfrac{1}{2}}} & & D_2
\ar[dl]^-{\ev_{{\tfrac{1}{2}}}} \\ & A \otimes B & }
\end{split}
\end{equation}
where 
$$D_1 = \{f \in C([0,\tfrac{1}{2}],A \otimes B) \mid f(0) \in A
\otimes B_1\}, \quad
D_2 = \{f \in C([\tfrac{1}{2},1],A \otimes B) \mid f(1) \in A \otimes
B_2\},$$
and where $\varphi_1$ and $\varphi_2$ are the restriction mappings.

We shall now use Proposition~\ref{prop:pull-back} to prove that $a \otimes 1$ 
belongs to the closure of the invertibles in $A \otimes E_{2,3}$. Condition
(iv) of Proposition~\ref{prop:pull-back} is satisfied by the assumption
on $a$ made in the first paragraph of the proof. The \Cs s $D_1$ and
$D_2$ are approximately divisible because $D_j \cong D_j \otimes
B_j$. It thus follows from  
Lemma~\ref{lm:K1} that condition (i) of Proposition~\ref{prop:pull-back}
is satisfied. The \Cs{} $A \otimes B$ is of stable rank one (cf.\
\cite{Ror:uhf}), whence $A \otimes B$ is strongly $K_1$-injective
(cf.\ the remarks below Lemma~\ref{lm:K1}). 

We have natural inclusions $A \otimes B_j \subseteq D_j$ (identifying 
an element in
$A \otimes B_j$ with a constant function), and the composition $A \otimes B_j
\to D_j \to A \otimes B$ is the inclusion mapping. Hence, to prove
that (iii) of Proposition~\ref{prop:pull-back} is satisfied, it suffices to
show that $K_1(A \otimes B)$ is generated by the images of the two mappings
$K_1(A \otimes B_j) \to K_1(A \otimes B)$, $j=1,2$. We have natural
identifications:
\begin{eqnarray*}
&K_1(A \otimes B_1) = K_1(A) \otimes \Z[1/2], \qquad 
K_1(A \otimes B_2) = K_1(A) \otimes \Z[1/3],& \\
& K_1(A \otimes B) = K_1(A) \otimes \Z[1/6]. &
\end{eqnarray*}
The desired identity now follows from the elementary fact that $\Z[1/2] +
\Z[1/3] = \Z[1/6]$.

Retaining the inclusion $A \otimes B_j \subseteq D_j$ from the previous 
paragraph, $\varphi_j(a \otimes 1_{E_{2,3}}) = a \otimes 1_{B_j}$.
Following \cite{Ror:uhf}, $a \otimes 1_{B_j}$ belongs to the closure of
the invertibles in $A \otimes B_j$ (and hence to the closure of the invertibles
in $D_j$) if (and only if) $\alpha_s(a)=0$; and $\alpha_s(a)=0$ for every
element $a$ in any unital, finite, simple \Cs{} $A$. It thus follows that
$\varphi_j(a)$ belongs to the closure of $\GL(D_j)$, $j=1,2$.
(If we had assumed that $A$ is stably finite, then we could have used
\cite[Corollary~6.6]{Ror:uhf} to conclude that the stable rank of 
$A \otimes B_j$ is one, which would have given us a more direct route
to the conclusion above.)
\end{proof}

\begin{theorem} \label{thm:sr1}
Every simple, unital, finite $\cZ$-absorbing \Cs{} 
has stable rank one.
\end{theorem}

\begin{proof}
Let $A$ be a simple, unital, finite \Cs{} such that $A$ is isomorphic
to $A \otimes \cZ$. Let $a \in A$ and $\ep>0$ be given. 
It follows from Lemma~\ref{lm:cts_field} that there is 
an invertible element $b \in A \otimes \cZ$ such that $\|a \otimes 1 -
b\| < \ep/2$. Let $\sigma_n \colon A \otimes \cZ \to A$ be as in
Lemma~\ref{lm:A_k} and choose $n$ such that $\|\sigma_n(a \otimes 1)-a\| <
\ep/2$. Then $\|a-\sigma_n(b)\| < \ep$, and $\sigma_n(b)$ is an
invertible element in $A$. 
\end{proof}

\section{The real rank of $\cZ$-absorbing \Cs s}

\noindent We conclude this paper with a result that describes when a
simple $\cZ$-absorbing \Cs{} is of real rank zero. A simple 
\emph{approximately divisible} 
\Cs{} is of real rank zero if and only if projections
separate quasitraces, as shown in \cite{BlaKumRor:apprdiv} (and, as
remarked earlier, each quasitrace on an exact \Cs{} is a trace by
\cite{Haa:quasi} and \cite{Kir:quasitraces}). It is
not true that any $\cZ$-absorbing \Cs{}, where quasitraces are being
separated by projections, is of real rank zero. The Jiang--Su
$\cZ$ itself is a counterexample. We must therefore require some
further properties, for 
example that the $K_0$-group is \emph{weakly divisible}: for each $g
\in K_0^+$ and for each $n \in \N$ there are $h_1,h_2 \in K_0^+$ such
that $g = nh_1+(n+1)h_2$. 

Let $A$ be a unital \Cs{}. Let $T(A)$ denote the simplex of
tracial states on $A$, and let $\Aff(T(A))$ denote the normed space of
real valued affine continuous functions on $T(A)$. Let $\rho \colon K_0(A) \to
\Aff(T(A))$ be the canonical map defined by $\rho(g)(\tau) =
K_0(\tau)(g)$.  

The result below is essentially contained in 
\cite[Theorem~7.2]{Ror:UHFII} (see also \cite{BlaKumRor:apprdiv}). 

\begin{proposition} \label{prop:rr0}
Let $A$ be an exact unital simple \Cs{} of stable rank one for which $W(A)$
is almost unperforated. Then $A$ is of real rank zero if and only if
$\rho(K_0(A))$ is uniformly dense in the normed space $\Aff(T(A))$.  
\end{proposition} 

\noindent The proof of \cite[Theorem~7.2]{Ror:UHFII} applies almost
verbatim. (At the point where we have a positive element $x \in A$ and
$\delta > 0$ such that $d_\tau(\langle f_{\delta/2}(x) \rangle) <
d_\tau(\langle f_{\delta/4}(x) \rangle)$ for all $\tau \in T(A)$,
then, because $\tau \mapsto
d_\tau(\langle y \rangle)$ defines an element in $\Aff(T(A))$ for
every $y \in A^+$, it follows by density of $\rho(K_0(A))$ in $\Aff(T(A))$
that there is an element $g \in K_0(A)$ such that $d_\tau(\langle
f_{\delta/2}(x) \rangle) < K_0(\tau)(g) < 
d_\tau(\langle f_{\delta/4}(x) \rangle)$ for all $\tau \in T(A)$. One
next uses weak unperforation of $K_0(A)$ (\cite{GongJiangSu:Z} or
Corollary~\ref{cor:K0-wup}) to conclude that $g$ is
positive, i.e., that $g = [q]$ for some projection $q$ in a matrix
algebra over $A$.)

\begin{theorem} \label{thm:RR0}  
The following conditions are equivalent for each unital, simple, exact,
finite, $\cZ$-absorbing \Cs{} $A$. 
\begin{enumerate}
\item $\RR(A) = 0$,
\item $\rho(K_0(A))$ is uniformly dense in $\Aff(T(A))$,
\item $K_0(A)$ is weakly divisible and projections in $A$ separate traces
  on $A$,
\end{enumerate}
\end{theorem}

\noindent
If $A$ is a simple, \emph{infinite}, $\cZ$-absorbing \Cs, then $A$ 
is purely infinite by \cite{GongJiangSu:Z}; and 
purely infinite \Cs s are of real rank zero by \cite{zhang:infsimp}. 

\begin{proof}
It follows from Theorem~\ref{thm:wup} that $W(A)$ is almost
unperforated, and from Theorem~\ref{thm:sr1} that the stable rank of
$A$ is one. We therefore get (ii) $\Rightarrow$ (i) from
Proposition~\ref{prop:rr0}. If (iii) holds, then for each element $g$
in $K_0(A)^+$ and for each natural number $n$ there exists $h \in
K_0(A)^+$ such that $nh \le g \le (n+1)h$. This implies that the uniform
closure of $\rho(K_0(A))$ in $\Aff(T(A))$ is a closed subspace which
separates points. Thus, by Kadison's Representation Theorem (see
\cite[II.1.8]{Alf:convex}), $\rho(K_0(A))$ is uniformly dense in
$\Aff(T(A))$, so (ii) holds. If (i) holds, then $A$ is weakly
divisible by \cite{PerRor:AF} and traces on $A$ are separated by
projections. 
\end{proof}

\noindent A \Cs{} is said to have property (SP) (``small
projections'') if each non-zero hereditary sub-\Cs{} contains a
non-zero projection.

\begin{corollary} \label{cor:rr0}
Let $A$ be a simple, unital, exact $\cZ$-absorbing \Cs{} with a 
\emph{unique} trace $\tau$. Then the following conditions are equivalent:
\begin{enumerate}
\item $\RR(A) = 0$,
\item $K_0(\tau)(K_0(A))$ is dense in $\R$,
\item $K_0(A)$ is weakly divisible,
\item $A$ has property (SP).
\end{enumerate}
\end{corollary}

\begin{proof}
The equivalence of (i), (ii), and (iii) follows immediately from
Theorem~\ref{thm:RR0}. The implication (i) $\Rightarrow$ (iv) is
trivial, and one easily sees that (iv) implies (ii). 
\end{proof}

\bibliographystyle{amsplain}

\begin{thebibliography}{10}

\bibitem{Alf:convex}
E.~Alfsen, \emph{{Compact convex sets and boundary integrals}},
  Springer-Verlag, New York, 1971, Ergebnisse der Mathematik und ihrer
  Grenzgebiete, Band 57.

\bibitem{BlaHan:quasitrace}
B.~Blackadar and D.~Handelman, \emph{{Dimension functions and traces on
  $C^*$-algebras}}, J. Funct. Anal. \textbf{45} (1982), 297--340.

\bibitem{BlaKumRor:apprdiv}
B.~Blackadar, A.~Kumjian, and M.~R{\o}rdam, \emph{{Approximately central matrix
  units and the structure of non-commutative tori}}, $K$-theory \textbf{6}
  (1992), 267--284.

\bibitem{BlaRor:extending}
B.~Blackadar and M.~R{\o}rdam, \emph{{Extending states on Preordered semigroups
  and the existence of quasitraces on $C^*$-algebras}}, J. Algebra \textbf{152}
  (1992), 240--247.

\bibitem{Cuntz:dimension}
J.~Cuntz, \emph{{Dimension functions on simple $C^*$-algebras}}, Math. Ann.
  \textbf{233} (1978), 145--153.

\bibitem{EilLorPed:semiproj}
S.~Eilers, T.~Loring, and G.~K. Pedersen, \emph{{Stability of anticommutation
  relations. An application of non-commutative CW complexes}}, J. Reine Angew.
  Math. \textbf{499} (1998), 101--143.

\bibitem{Ell:torsion}
G.~A. Elliott, \emph{{Dimension groups with torsion}}, Int. J. Math. \textbf{1}
  (1990), 361--380.

\bibitem{Ell:classprob}
\bysame, \emph{{The classification problem for amenable $C^*$-algebras}},
  Proceedings of the International Congress of Mathematicians (Basel), vol.
  1,2, Birkh\"auser, 1995, pp.~922--932.

\bibitem{GongJiangSu:Z}
G.~Gong, X.~Jiang, and H.~Su, \emph{{Obstructions to $\mathcal{Z}$-stability
  for unital simple $C^*$-algebras}}, Canadian Math. Bull. \textbf{43} (2000),
  no.~4, 418--426.

\bibitem{GooHan:extending}
K.~R. Goodearl and D.~Handelman, \emph{{Rank functions and $K_0$ of regular
  rings}}, J. Pure Appl. Algebra \textbf{7} (1976), 195--216.

\bibitem{Haa:quasi}
U.~Haagerup, \emph{{Every quasi-trace on an exact $C^*$-algebra is a trace}},
  preprint, 1991.

\bibitem{JiaSu:Z}
X.~Jiang and H.~Su, \emph{{On a simple unital projectionless $C^*$-algebra}},
  American J. Math. \textbf{121} (1999), no.~2, 359--413.

\bibitem{Kir:spi}
E.~Kirchberg, \emph{{Permanence properties of purely infinite $C^*$-algebras}},
  in preparation.

\bibitem{Kir:quasitraces}
\bysame, \emph{{On the existence of traces on exact stably projectionless
  simple $C^*$-algebras}}, Operator Algebras and their Applications (P.~A.
  Fillmore and J.~A. Mingo, eds.), Fields Institute Communications, vol.~13,
  Amer. Math. Soc., 1995, pp.~171--172.

\bibitem{KirRor:pi}
E.~Kirchberg and M.~R{\o}rdam, \emph{{Non-simple purely infinite
  $C^*$-algebras}}, American J. Math. \textbf{122} (2000), 637--666.

\bibitem{KirRor:pi2}
\bysame, \emph{{Infinite non-simple $C^*$-algebras: absorbing the Cuntz algebra
  $\mathcal{O}_\infty$}}, Advances in Math. \textbf{167} (2002), no.~2,
  195--264.

\bibitem{Lin:amenable}
H.~Lin, \emph{{An introduction to the classification of amenable
  {$C^*$}-algebras}}, World Scientific Publishing Co. Inc., River Edge, NJ,
  2001.

\bibitem{PerRor:AF}
F.~Perera and M.~R{\o}rdam, \emph{{AF-embeddings into $C^*$-algebras of real
  rank zero}}, to appear in J. Funct. Anal.

\bibitem{Rfl:sr}
M.~Rieffel, \emph{{Dimension and stable rank in the $K$-theory of
  $C^*$-algebras}}, Proc. London Math. Soc. \textbf{46} (1983), no.~(3),
  301--333.

\bibitem{Ror:unitary}
M.~R{\o}rdam, \emph{{Advances in the theory of unitary rank and regular
  approximation}}, Ann. of Math. \textbf{128} (1988), 153--172.

\bibitem{Ror:uhf}
\bysame, \emph{{On the Structure of Simple $C^*$-algebras Tensored with a
  UHF-Algebra}}, J. Funct. Anal. \textbf{100} (1991), 1--17.

\bibitem{Ror:UHFII}
\bysame, \emph{{On the Structure of Simple $C^*$-algebras Tensored with a
  UHF-Algebra, II}}, J. Funct. Anal. \textbf{107} (1992), 255--269.

\bibitem{Ror:encyc}
\bysame, \emph{{Classification of Nuclear, Simple $C^*$-algebras}},
  {Classification of Nuclear $C^*$-Algebras. Entropy in Operator Algebras}
  (J.~Cuntz and V.~Jones, eds.), vol. 126, Encyclopaedia of Mathematical
  Sciences. Subseries: Operator Algebras and Non-commutative Geometry, no. VII,
  Springer Verlag, Berlin, Heidelberg, 2001, pp.~1--145.

\bibitem{Ror:simple}
\bysame, \emph{{A simple $C^*$-algebra with a finite and an infinite
  projection}}, Acta Math. \textbf{191} (2003), 109--142.

\bibitem{Toms:example}
A.~Toms, \emph{{On the Independence of $K$-theory and Stable Rank for Simple
  $C^*$-algebras}}, preprint.

\bibitem{TomsWin:Z}
A.~Toms and W.~Winter, in preparation.

\bibitem{Vil:perforation}
J.~Villadsen, \emph{{Simple $C^*$-algebras with perforation}}, J. Funct. Anal.
  \textbf{154} (1998), no.~1, 110--116.

\bibitem{zhang:infsimp}
S.~Zhang, \emph{{A property of purely infinite simple $C^*$-algebras}}, Proc.
  Amer. Math. Soc. \textbf{109} (1990), 717--720.

\end{thebibliography}

\providecommand{\bysame}{\leavevmode\hbox to3em{\hrulefill}\thinspace}
\providecommand{\MR}{\relax\ifhmode\unskip\space\fi MR }
\providecommand{\MRhref}[2]{%
  \href{http://www.ams.org/mathscinet-getitem?mr=#1}{#2}
}
\providecommand{\href}[2]{#2}

\vspace{.5cm}

\noindent{\sc Department of Mathematics, University of Southern
  Denmark, Odense,
  Campusvej~55, 5230 Odense M, Denmark}

\vspace{.3cm}

\noindent{\sl E-mail address:} {\tt mikael@imada.sdu.dk}\\
\noindent{\sl Internet home page:}
{\tt www.imada.sdu.dk/$\,\widetilde{\;}$mikael/welcome} \\

\end{document}